\documentclass{article}
\usepackage{graphicx} 
\usepackage{amsmath,amssymb}
\usepackage{amsthm}
\usepackage{amsfonts}
\usepackage{txfonts}
\usepackage{setspace}
\usepackage{xcolor}
\usepackage{appendix}
\usepackage[utf8]{inputenc}
\usepackage{indentfirst}
\usepackage{listings}
\usepackage{url}

\newcommand{\ZZ}{\mathbb Z}

\newcommand{\Qc}{{\mathcal Q}}
\newcommand{\id}{{\mathrm {id}}}
\newcommand{\auto}{{\mathrm {Aut}}}
\newcommand{\inn}{{\mathrm {Inn}}}

\newcommand{\fix}{{\mathrm {Fix}}}
\newcommand{\gaq}{{\mathrm {GAQ}}}

\theoremstyle{plain}
\newtheorem{thm}{Theorem}[section]
\newtheorem{prop}{Proposition}[section]
\newtheorem{prob}{Problem}[section]
\theoremstyle{definition}
\newtheorem{defi}{Definition}[section]

\theoremstyle{remark}
\newtheorem{rem}{Remark}[section]

\numberwithin{equation}{section}

\makeatletter
\newcommand{\subjclass}[2][1991]{%
  \let\@oldtitle\@title%
  \gdef\@title{\@oldtitle\footnotetext{#1 Mathematics subject classification: #2}}%
}
\newcommand{\keywords}[1]{%
  \let\@@oldtitle\@title%
  \gdef\@title{\@@oldtitle\footnotetext{Key words and phrases: #1.}}%
}
\makeatother

\title{Classification of generalized Alexander quandles}
\author{Akihiro Higashitani, Seiichi Kamada, Jin Kosaka and Hirotake Kurihara}
\date{}

\begin{document}

\maketitle

\begin{abstract}
The aim of this paper is to provide a new characterization of isomorphism classes of 
generalized Alexander quandles in terms of the underlying groups and their automorphisms. 
This extends the previous result \cite[Theorem 1.4]{HK2022}. Additionally, we compute the number of generalized Alexander quandles up to quandle isomorphism arising from groups up to order $127$ and their group automorphisms. 
\end{abstract}

\section{Introduction}

\subsection{Motivation}

Quandles \cite{J1982, M1982} are algebraic systems equipped with three axioms that correspond to Reidemeister moves in knot theory. 
According to Joyce, any quandle can be expressed by a disjoint union of homogeneous quandles (\cite[Theorem 7.2]{J1982}). 
Here, we say that a quandle $Q$ is \textit{homogeneous} if the automorphism group of $Q$ acts transitively on $Q$. 
Furthermore, it has been shown that any homogeneous quandle can be constructed by a so-called \textit{quandle triplet} $(G,~H,~\psi)$ (\cite[Theorem 7.1]{J1982}). 

The primary focus of this paper is a special kind of homogeneous quandles called \textit{generalized Alexander quandles}, which are quandles arising from groups.  
Alexander quandles, which are generalized Alexander quandles arising from abelian groups, have been extensively studied in the context of knot theory. 
Finite Alexander quandles have been classified up to quandle isomorphism by Nelson \cite{N2003}. 

Recently, the first named and fourth named authors established a kind of characterization theorem of finite generalized Alexander quandles of finite groups under certain assumptions $(P1)$ and $(P2)$ (detailed in Section~\ref{sec:HK}.) 
The main objective of this paper is to extend the result  of \cite{HK2022}. 

\subsection{Main results}

The following is the main theorem of this paper: 
\begin{thm}[~Main Result~]\label{thm:main}
    Let $G$ and $G'$ be groups with their identity elements $e$ and $e'$, respectively. 
    Let $\psi \in \auto(G),~\psi' \in \auto(G')$ and let
     \[Q = (Q(G,\psi),*),~Q' = (Q(G',\psi'),*'),~P = P(Q)\text{ and }P' = P(Q').\]
     Then the following conditions $(1)$ and $(2)$ are equivalent: 
    \begin{itemize}
        \item[$(1)$] $Q$ and $Q'$ are isomorphic as quandles.  
        \item[$(2)$] There exists a group isomorphism $h : P \rightarrow P'$ such that the following two conditions $(A)$ and $(B)$ are satisfied: 
  \begin{itemize}
        \item [$(A)$] $h \circ \psi|_P = \psi'|_{P'} \circ h$; 
        \item [$(B)$] there exist representatives $A\subset G$ and $A'\subset G'$ of $G/P$ and $G'/P'$, respectively, and a bijection $k : A \rightarrow A'$ between $A$ and $A'$ such that $h(e*a) = e'*' k(a)$ holds for any $a \in A$. 
  \end{itemize}  
  \end{itemize}
\end{thm}
The notation used in this paper will be explained in Sections~\ref{sec:pre} and \ref{sec:HK}.
The proof of Theorem~\ref{thm:main} will be given in Section~\ref{sec:proof}. 
\begin{rem}
    The statement of Theorem~\ref{thm:main} does not assume the finiteness of groups. 
\end{rem}

By \cite{HK2022}, generalized Alexander quandles of finite groups up to order $15$ are completely classified. 
An observation of those of order $16$ is also made, but the main result of \cite{HK2022} is not applicable for two cases. 
By Theorem~\ref{thm:main}, we can completely classify not only all generalized Alexander quandles of finite groups of order $16$, but also those with orders at most $127$ (with some additional exceptional treatments). 
See Section~\ref{sec:app}. 

\subsection*{Acknowledgements}

The authors are grateful to S.~Nelson and H.~Tamaru for useful comments on an earlier draft of the paper.
The first named author is partially supported by JSPS KAKENHI Grant Number JP20K03513. 
The second named author is supported by JP19H01788 and JP23H05437. 
The fourth named author is supported by JP20K03623.

\bigskip

\section{Preliminaries}\label{sec:pre}

In this section, we prepare notation and recall notions on quandles and groups used throughout the paper. 

\subsection{Definition of quandles}

First, let us recall a definition of quandle. 
A non-empty set $Q$ equipped with a binary operation $*$ satisfying the following three axioms is called a \textit{quandle}: 
\begin{enumerate}
    \item[(Q1)] $x*x=x$ for any $x \in Q$; 
    \item[(Q2)] for any $x \in Q$, the map defined by $s_x : Q \ni y \mapsto y * x \in Q$ is bijective; 
    \item[(Q3)] for any $x,y,z \in Q$, we have $(x*y)*z=(x*z)*(y*z)$. 
\end{enumerate}
The map appearing in (Q2) is called a \textit{point symmetry}. 

\subsection{Homogeneous quandles and generalized Alexander quandles}

Let $G$ be a group and let $\psi$ be an automorphism of $G$. 
Moreover, we take a subgroup $H$ of $G$ satisfying $\psi(h)=h$ for any $h \in H$. 
We define a binary operation $*$ on the set of left cosets $G/H$ defined by $xH*yH:=y \psi(y^{-1}x)H$. 
Then this binary operation becomes well-defined and $(G/H,*)$ gives us a quandle. 
The triplet $(G,H,\psi)$ providing a new quandle $(G/H,*)$ is called a \textit{quandle triplet}, and we denote this quandle by $Q(G,H,\psi)$. 
We notice that for any group $G$ with its identity element $e$ and any automorphism $\psi$ of $G$, 
the triplet $(G,\{e\},\psi)$ always becomes a quandle triplet. 

\begin{defi}
Let $G$ be a group and let $\psi$ be its automorphism. 
Considering a binary operation $*$ on $G$ defined by $x*y:=y\psi(y^{-1}x)$, we obtain a quandle $(G,*)$, denoted by $Q(G,\psi)$ or $(Q(G,\psi), *)$. 
We call this quandle $Q(G,\psi)$ a \textit{generalized Alexander quandle} (or a \textit{principal quandle}, see~\cite{Bonatto2020}). 
\end{defi}

In the case where $G$ is abelian, we call $Q(G,\psi)$ an \textit{Alexander quandle}, 
which has been studied especially in the context of knot theory. 
Generalized Alexander quandles can be regarded as a generalization of it. 

\subsection{Notation and terminologies on quandles and groups}

In what follows, we use the notation $(Q,*)$ for quandles if we would like to emphasize the equipped binary operation; otherwise, we simply write $Q$. 

Let $(Q,*)$ and $(Q',*')$ be quandles. 
A map $f:Q \rightarrow Q'$ is said to be a \textit{quandle homomorphism} if $f(x*y)=f(x)*'f(y)$ holds for any $x,y \in Q$. 
In addition to this, if $f$ is bijective, then we call $f$ a \textit{quandle isomorphism}. 
If there exists a quandle isomorphism between quandles $Q$ and $Q'$, then we say that $Q$ and $Q'$ are isomorphic as quandles, which we denote by $Q \cong Q'$. 
\begin{itemize}
    \item Let $\auto(Q)$ be the group consisting of all quandle automorphisms on $Q$, which we call the \textit{automorphism group} of $Q$. 
    Note that the point symmetry $s_x$ for each $x \in Q$ belongs to $\auto(Q)$. 
    \item Let $\inn(Q)$ be the subgroup of $\auto(Q)$ generated by $\{s_x : x \in Q\}$, which we call the \textit{inner automorphism group} of $Q$. 
    More explicitly, we have $\inn(Q)=\{s_{x_n}^{\pm 1} \circ s_{x_{n-1}}^{\pm 1} \circ \cdots \circ s_{x_1}^{\pm 1} : n \in \ZZ_{> 0}, x_1,\ldots,x_n \in Q\}$. 
    In the case where $Q$ is finite, we have $\inn(Q)=\{\id_Q, s_{x_n} \circ s_{x_{n-1}} \circ \cdots \circ s_{x_1} : n \in \ZZ_{> 0}, x_1,\ldots,x_n \in Q\}$. 
\end{itemize}

We say that a quandle $Q$ is \textit{homogeneous} 
if the automorphism group $\auto(Q)$ of $Q$ acts on $Q$ transitively. 
According to \cite[Section 7]{J1982}, the following one-to-one correspondence is known: 
\begin{align*}
    \{\text{homogeneous quandles}\}/ \cong \;\;\overset{1 : 1}{\longleftrightarrow} \;\; 
    \{Q(G,H,\psi) : (G,H,\psi) \text{ is a quandle triplet}\}/ \cong, 
\end{align*}
where $/\cong$ stands for that we consider this correspondence up to quandle isomorphism. 

\begin{rem}\label{rem:e}
    Let $Q$ and $Q'$ be homogeneous quandles which are isomorphic. 
    Then, for any $x \in Q$ and $x' \in Q'$, we can always find a quandle isomorphism $f: Q \rightarrow Q'$ with $f(x)=x'$. 
    In fact, even if $f(x)=x'' \neq x'$, since $\auto(Q')$ acts on $Q'$ transitively, 
    we can find $g \in \auto(Q')$ with $g(x'')=x'$, so we may take an isomorphism $g \circ f$ between $Q$ and $Q'$ instead of $f$. 

    Therefore, for any generalized Alexander quandles $Q(G,\psi)$ and $Q(G',\psi')$ which are isomorphic, 
    we may assume the existence of a quandle isomorphism $f:Q(G,\psi) \rightarrow Q(G',\psi')$ with $f(e)=e'$, 
    where $G$ and $G'$ are groups with their identity elements $e$ and $e'$, respectively.  
\end{rem}

\medskip

We also collect some terminologies and notation concerning with groups. Let $G$ be a group. 
\begin{itemize}
    \item We denote the automorphism group of $G$ by $\auto(G)$. Remark that this notation is the same as that for quandles. 
    \item Given $\psi \in \auto(G)$, let $\fix(\psi,G)=\{g \in G : \psi(g)=g\}$. 
    \item Let $N,H$ be groups and let $\sigma : H \rightarrow \auto(N)$, $h \mapsto \sigma_h$ be a group homomorphism. 
    The (outer) \textit{semidirect product} $N \rtimes_\sigma H$ of $N$ and $H$ with respect to $\sigma$ is a group 
    which is the direct product $N \times H$ as a set equipped with the operation defined by $(n_1,h_1)(n_2,h_2)=(n_1\sigma_{h_1}(n_2),h_1h_2)$. 
\end{itemize}

\subsection{Problem}

Given a group $G$, let $\Qc(G)$ denote the set of all generalized Alexander quandles $Q(G,\psi)$ with $\psi \in \auto(G)$ up to quandle isomorphism. Namely, 
\[
\Qc(G)=\{Q(G,\psi) : \psi \in \auto(G)\}/\cong. 
\]

\begin{prob}[ {\cite[Problem 3.5]{HK2020}} ]\label{prob:main}
    Determine $\Qc(G)$ for a given group $G$. 
\end{prob}
This problem has some certain partial answers in terms of the automorphism group of groups. 
For example, \cite[Corollary 2.2]{N2003} claims that $\Qc(G)$ can be characterized in terms of just numerical information in the case where $G$ is finite and abelian.
Moreover, it is proved in \cite[Theorem 1.2]{HK2022} that if $G$ is a finite simple group, then 
$\Qc(G)$ one-to-one corresponds to the conjugacy classes of $\auto(G)$. 
Using this, we can also get this one-to-one correspondence in the case of symmetric groups (\cite[Corollary 1.3]{HK2022}), though it is not simple. 
Theorem~\ref{thm:main} asserts a new kind of complete solution of Problem~\ref{prob:main}. 

\section{Previous results on generalized Alexander quandles}\label{sec:HK}

In this section, we recall some results on generalized Alexander quandles, especially from \cite{HK2022}. 

Let $G$ be a group, $\psi$ an automorphism of $G$, and $Q=Q(G,\psi)$ the generalized Alexander quandle. 
We define $P=P(Q)$ as the orbit of the identity element $e$ in $Q$ under the action of $\inn(Q)$, i.e., 
\[
P=P(Q)=
\{s_{x_n}^{\pm 1} \circ s_{x_{n-1}}^{\pm 1} \circ \cdots \circ s_{x_1}^{\pm 1}(e) : n \in \ZZ_{> 0}, x_1,\ldots,x_n \in Q\}.
\]
Then $P$ enjoys a special property. In particular, the following statements hold: 

\begin{prop}[ {\cite[Propositions 3.2 and 3.3]{HK2022}} ]\label{prop:normal}
The subset $P$ of $G$ is a subquandle of $Q$ and a normal subgroup of $G$. 
\end{prop}

\begin{thm}[ {\cite[Theorem 3.10]{HK2022}, \cite[Lemma 2.4]{T23}} ]\label{thm:HK}
Let $G$ and $G'$ be groups with their identity elements $e$ and $e'$, respectively. 
Let $\psi \in \mathrm{Aut}(G),~\psi' \in \mathrm{Aut}(G')$ and let 
\[Q = Q(G,\psi),~Q' = Q(G',\psi'),~P = P(Q)\text{ and }P' = P(Q').\]
Assume that there is a quandle isomorphism $f:Q \rightarrow Q'$ with $f(e)=e'$. 
Then the following hold: 
\begin{itemize}
    \item[(i)] $f \circ \psi = \psi' \circ f$; 
    \item[(ii)] $f|_P : P \rightarrow P'$ is a group isomorphism; 
    \item[(iii)] $f|_P : Q(P,\psi|_P) \rightarrow Q(P',\psi'|_{P'})$ is a quandle isomorphism; 
    \item[(iv)] $f(xP)=f(x)P'$ holds for any $x \in G$. 
\end{itemize}
\end{thm}

\begin{rem}
    In \cite[Theorem 3.10]{HK2022}, the finiteness of groups is assumed, but Proposition~\ref{prop:normal} and Theorem~\ref{thm:HK} do not need this assumption. 
    In fact, even in the case of infinite groups, we see that 
    Proposition~\ref{prop:normal} and the statements of (i), (iii) and (iv) of Theorem~\ref{thm:HK} hold by almost the same proofs with those given in \cite{HK2022}. 
    Moreover, the statement (ii) for (non-necessarily finite) groups is nothing but \cite[Lemma 2.4]{T23}. 
\end{rem}

\bigskip

The following proposition is useful for study of generalized Alexander  quandles up to isomorphism. It asserts that 
a conjugacy class of $\auto(G)$ determines an element of $\Qc(G)$. 

\begin{prop}[ {\cite[Proposition 3.5]{HK2022}} ]
    Let $\psi,\psi' \in \auto(G)$ and let $Q=Q(G,\psi)$ and $Q'=Q(G,\psi')$. 
    If $\psi$ and $\psi'$ are conjugate in $\auto(G)$, then $Q \cong Q'$ holds. 
\end{prop}

\medskip

Now, we recall the main result of \cite{HK2022}. For this, we have to introduce more notions. 
Given a group $G$ and its automorphism $\psi$, let $Q=Q(G,\psi)$ and $P=P(Q)$. We can consider a generalized Alexander quandle  $Q(P,\psi|_P)$. 
Let $P^2$ be $P(Q(P,\psi|_P))$.  
By Proposition~\ref{prop:normal}, $P^2$ is a normal subgroup of $P$. 

We consider the following two properties (P1) and (P2). 
\begin{itemize}
\item[(P1)] $P^2$ is a normal subgroup of $G$. 
\item[(P2)] $P^2=\{s_p(e) : p \in P\}$ holds. 
\end{itemize}

These two properties are independent each other and become invariants on generalized Alexander quandles (see \cite[Proposition 3.16]{HK2022}). 

\begin{thm}[ {\cite[Theorem 1.4]{HK2022}} ]\label{thm:HK2}
    Let $G$ and $G'$ be finite groups, $\psi \in \mathrm{Aut}(G)$, $\psi' \in \mathrm{Aut}(G')$ and let 
    \[Q = (Q(G,\psi),*),~Q' = (Q(G',\psi'),*'),~P = P(Q) \text{ and }P' = P(Q').\]
    Assume that $Q$ and $Q'$ satisfy (P1) and (P2). Then the following conditions (1) and (2) are equivalent: 
    \begin{itemize}
        \item[$(1)$] $Q \cong Q'$. 
        \item[$(2)$] The following three conditions are satisfied: 
    \begin{itemize}
        \item[(i)] $|G|=|G'|$; 
        \item[(ii)] $|\fix(\psi,G)|=|\fix(\psi',G')|$; 
        \item[(iii)] there exists a group isomorphism $h: P \rightarrow P'$ such that the following $(A)$ and $(B)$ are satisfied: 
        \begin{itemize}
            \item[$(A)$] $h \circ \psi|_P = \psi'|_{P'} \circ h$; 
            \item[$(B)$] for any $a \in G$, there exists $a' \in G'$ such that $h(e*a)=e'*' a'$.  
        \end{itemize}
    \end{itemize}  
    \end{itemize}
\end{thm}

Note that the finiteness of the groups $G$ and $G'$ is necessary in this theorem since we have to take the orders of the groups into account, 
while Theorem~\ref{thm:main} does require neither the finiteness condition on $G$ and $G'$ nor the conditions (P1) and (P2) on $Q$ and $Q'$.

\section{Proof of Theorem~\ref{thm:main}}\label{sec:proof}

This section is devoted to proving Theorem~\ref{thm:main}. 

\bigskip

\noindent
(1) $\Rightarrow$ (2): 

Let $f$ be a quandle isomorphism $f:Q \rightarrow Q'$ with $f(e)=e'$. (See Remark~\ref{rem:e}.) 
Let $h$ be the restriction of $f$ into $P$, i.e., $h=f|_P: P \rightarrow P'$. 
Then, by Theorem~\ref{thm:HK} (i) and (ii), we see that the condition (A) is satisfied. 

Take an arbitrary representative $A\subset G$ of $G/P$, let $A'=f(A)$ 
and let $k$ be the restriction of $f$ into $A$, i.e., $k=f|_A:A \rightarrow A'$. 
Then it follows from the bijectivity of $f$ (i.e., $k$) together with Theorem~\ref{thm:HK} (iv) that $A'$ also becomes a representative of $G'/P'$. 
Since $e*a$ is an element of $P$, we see that $h(e*a)=f(e*a)=f(e)*' f(a)=e'*' k(a)$ holds for any $a \in A$. 
Hence, (B) is also satisfied. 

\medskip

\noindent
(2) $\Rightarrow$ (1): 

Let $h:P \rightarrow P'$ and $k:A \rightarrow A'$ be a group isomorphism and a bijection between representatives $A$ and $A'$, respectively, 
satisfying the conditions (A) and (B). 
Given any $x \in G$, we can uniquely choose $a_x \in A$ and $p_x \in P$ such that $x=a_xp_x$. 
Since $P$ is a normal subgroup of $G$ (see Proposition~\ref{prop:normal}), we know that $a_xp_aa_x^{-1} \in P$ holds. 

Now, we define a map $f: G \rightarrow G'$ by setting 
\[
f(x):=h(a_xp_xa_x^{-1})k(a_x) \in G' 
\]
for any $x=a_xp_x \in G$. 
The bijectivity of $f$ follows from that of $h$ and $k$. 
In what follows, we prove that $f$ is a quandle homomorphism. 

Since $h(a_xp_xa_x^{-1}) \in P'$ and $P'$ is a normal subgroup of $G'$,
$k(a_x)^{-1}h(a_xp_xa_x^{-1})k(a_x)$ belongs to $P'$.
Hence the coset decomposition of $f(x)=a'_{f(x)} p'_{f(x)}$ becomes $a'_{f(x)}=k(a_x)$ and $p'_{f(x)}=k(a_x)^{-1}h(a_xp_xa_x^{-1})k(a_x)$.

For any $x =a_xp_x$ and $y=a_yp_y$ in $G$, 
we have
\[
x*y= a_yp_y\psi((a_yp_y)^{-1}a_xp_x)
    = a_x \cdot a_x^{-1} a_y p_y\psi((a_x^{-1} a_y p_y)^{-1}) \psi(p_x).
\]
Since we have $a_x^{-1} a_y p_y\psi((a_x^{-1} a_y p_y)^{-1}) = 
e*(a_x^{-1} a_y p_y) 
\in P$ and $\psi(p_x)\in P$, 
the coset decomposition of $x*y=a_{x*y} p_{x*y}$ becomes 
\[
a_{x*y}=a_x
\]
and 
\[p_{x*y}=a_x^{-1} a_y p_y\psi((a_x^{-1} a_y p_y)^{-1}) \psi(p_x).
\]
Thus, we have the following:
\begin{align*}
    f(x*y)
    &=h(a_{x*y}p_{x*y}a_{x*y}^{-1})k(a_{x*y})\\
    &=h(a_x \cdot a_x^{-1} a_y p_y\psi((a_x^{-1} a_y p_y)^{-1}) \psi(p_x) \cdot a_x^{-1}) k(a_x)\\
    &= h(a_y p_y~\psi(p_y^{-1} a_y^{-1})~\psi(a_x p_x) a_x^{-1})~k(a_x) \\
    &= h(a_y p_y a_y^{-1} \cdot a_y~\psi(p_y^{-1} a_y^{-1})~(a_x~\psi(p_x^{-1} a_x^{-1}))^{-1})~k(a_x) \\
    &= h(a_y p_y a_y^{-1})~h(a_y~\psi(p_y^{-1} a_y^{-1}))~h(a_x~\psi(p_x^{-1} a_x^{-1}))^{-1}~k(a_x) \\
    &= h(a_y p_y a_y^{-1})~h(a_y~\psi(a_y^{-1}) \cdot \psi(a_y p_y^{-1} a_y^{-1}))~h(\psi(a_x p_x a_x^{-1}) \cdot \psi(a_x)~a_x^{-1})~k(a_x) \\
    &= h(a_y p_y a_y^{-1})~h(a_y~\psi(a_y^{-1}))~h \circ \psi(a_y p_y^{-1} a_y^{-1} a_x p_x a_x^{-1})~h(\psi(a_x) a_x^{-1})~k(a_x).
\end{align*}

On the other hand, for the coset decomposition of 
$ f(x) *' f(y) = a'_{f(x) *' f(y)} p_{f(x) *' f(y)}'$, 
we also have $a'_{f(x) *' f(y)} =a'_{f(x)}=k(a_x)$ and
\begin{align*}
    p_{f(x) *' f(y)}'
    &={a_{f(x)}'}^{-1} {a_{f(y)}}^\prime {p_{f(y)}}^\prime \psi' (({a_{f(x)}^\prime}^{-1} a_{f(y)}'p_{f(y)}')^{-1}) \psi'(p_{f(x}') \\
    &=k(a_x)^{-1} k(a_y) \cdot k(a_y)^{-1}h(a_yp_ya_y^{-1})k(a_y)
    \psi'((k(a_x)^{-1} k(a_y) \cdot k(a_y)^{-1}h(a_yp_ya_y^{-1})k(a_y))^{-1})\\
    &\quad\quad \cdot \psi'(k(a_x)^{-1}h(a_xp_xa_x^{-1})k(a_x))\\
    &=k(a_x)^{-1} h(a_yp_ya_y^{-1})k(a_y)
    \psi'(k(a_y)^{-1} h(a_y p_y^{-1}a_y^{-1})h(a_x p_x a_x^{-1})k(a_x)).
\end{align*}
Therefore, we obtain
\begin{align*}
    f(x) *' f(y) 
    &= a'_{f(x) *' f(y)} p'_{f(x) *' f(y)}\\
    &=h(a_yp_ya_y^{-1})k(a_y) \psi'(k(a_y)^{-1} h(a_yp_y^{-1}a_y^{-1})h(a_xp_xa_x^{-1})k(a_x))\\
    &= h(a_y p_y a_y^{-1})~k(a_y)~\psi'(k(a_y)^{-1})~\psi'(h(a_y p_y^{-1} a_y^{-1} a_x p_x a_x^{-1}))~\psi'(k(a_x)) \\
    &= h(a_y p_y a_y^{-1})~k(a_y)~\psi'(k(a_y)^{-1})~h \circ \psi(a_y p_y^{-1} a_y^{-1} a_x p_x a_x^{-1})~\psi'(k(a_x))\;\; 
    (\text{by (A)}).
\end{align*}
We consider $k(a_y)~\psi'(k(a_y)^{-1})$ and $\psi'(k(a_x))$ in the above equation.
By (B), we have
\[k(a_y)~\psi'(k(a_y)^{-1})=e' *' k(a_y)=h(e * a_y)=h(a_y \psi (a_y^{-1}))\]
and
\begin{align*}
    \psi'(k(a_x))
    &=
    (k(a_x) \psi'(k(a_x)^{-1}))^{-1} ~ k(a_x) = (e'*' k(a_x))^{-1} ~ k(a_x) = h(e* a_x)^{-1} ~ k(a_x) \\
    &= h(a_x \psi(a_x^{-1}))^{-1} ~ k(a_x) = h( \psi(a_x) a_x^{-1}) ~ k(a_x).
\end{align*}
Thus, we have
\[f(x) *' f(y) =
h(a_y p_y a_y^{-1})~h(a_y~\psi(a_y^{-1}))~h \circ \psi(a_y p_y^{-1} a_y^{-1} a_x p_x a_x^{-1})~h(\psi(a_x) a_x^{-1})~k(a_x)
=f(x*y).\]

Therefore, we conclude that $f$ is a quandle isomorphism, as required. 

\bigskip

\section{Classification of generalized Alexander quandles arising from finite groups with small order}\label{sec:app}

In this section, we apply Theorem~\ref{thm:main} to obtain the complete classification of 
generalized Alexander quandles arising from finite groups with small orders. 

Given a positive integer $n$, let 
\begin{align*}
\Qc_\gaq(n):&=\{Q(G,\psi) : G \text{ is a group of order }n, \; \psi \in \auto(G)\}/\cong. 
\end{align*}
What we will work in this section is precisely to determine $|\Qc_\gaq(n)|$ for $n \leq 127$. 

\subsection{The case $n \leq 15$}
In \cite[Section 6]{HK2022}, classification of generalized Alexander quandles arising from finite groups $G$ is discussed. 
For example, the cases where $|G|=p,2p$ or $p^2$ are investigated theoretically, where $p$ is a prime. 
Moreover, it follows from \cite[Corollary 1.7]{HK2022} that any generalized Alexander quandles arising from cyclic groups of order $2n$ 
can be realized as those arising from dihedral groups. 
By additional case-by-case discussions, we can eventually get the complete classification of $\Qc_\gaq(n)$ for $n \leq 15$. 

\subsection{The case $n=16$}
Most parts of $\Qc_\gaq(16)$ are also discussed in \cite[Section 6]{HK2022}, 
although two cases $\{Q_1,Q_1'\}$ and $\{Q_2,Q_2'\}$ are remained to be open. 
In this subsection, we figure these cases out by using Theorem~\ref{thm:main}. 

\bigskip

Let $Q_8$ be the quaternion group 
$\langle \overline{e}, \mathbf{i},\mathbf{j},\mathbf{k} : 
\overline{e}^2 = e, 
\mathbf{i}^2 = \mathbf{j}^2 = \mathbf{k}^2 = \mathbf{i}\mathbf{j}\mathbf{k} = \overline{e}
\rangle$.  
We denote the elements 
$e$, $\mathbf{i}$, $\mathbf{j}$, $\mathbf{k}$, 
$\overline{e}$, $\overline{e}\mathbf{i}$, $\overline{e}\mathbf{j}$,  
$\overline{e}\mathbf{k}$ by 
$1, \mathbf{i},\mathbf{j},\mathbf{k}, -1, -\mathbf{i}, -\mathbf{j}, -\mathbf{k}$, respectively, so that 
$Q_8 = 
\{ 1, \mathbf{i},\mathbf{j},\mathbf{k}, -1, -\mathbf{i}, -\mathbf{j}, -\mathbf{k}
\}$. 

Let $\ZZ/2\ZZ=\{\bar{0},\bar{1}\}$ be the group of order $2$. 
Let $\sigma:\ZZ/2\ZZ \rightarrow \auto(Q_8)$ be defined by $i \mapsto \sigma_i$, where 
$\sigma_{\bar{0}} = \id_{Q_8}$ and $\sigma_{\bar{1}}:(\mathbf{i},\mathbf{j},\mathbf{k}) \mapsto (\mathbf{i},-\mathbf{j},-\mathbf{k})$. 
Then $\sigma_i$ can be written as the inner automorphism $\sigma_i(a)=\mathbf{i}^i a \mathbf{i}^{-i}$ for $a \in Q_8$. 
Note that $\mathbf{i}^i a \mathbf{i}^{-i}$ is well-defined on $\ZZ/2\ZZ$.
We define the groups $G$ and $G'$ as follows: 
\[G=Q_8 \times \ZZ/2\ZZ \;\;\text{ and }\;\;G'=Q_8 \rtimes_\sigma \ZZ/2\ZZ. \]
Let $\psi$ be the group automorphism of $Q_8$ 
defined by $\psi(\mathbf{i})=\mathbf{j}$ and $\psi(\mathbf{j})=\mathbf{k}$,
i.e., $\psi$ cyclically permutes $\mathbf{i}, \mathbf{j}, \mathbf{k}$.
Moreover, we take the following maps $\psi_1$ on $G$ and $\psi_1'$ on $G'$, respectively:
\begin{align*}
    \psi_1 (a,i) =(\psi(a),i) \;\;\text{and}\;\; \psi_1' (a,i) =(\psi(a \mathbf{i}^i) \mathbf{i}^{-i},i). 
\end{align*}
Note that we have 
\[\psi(a \mathbf{i}^i) \mathbf{i}^{-i}=\begin{cases}\psi(a) &\text{ if }i=\bar{0}, \\ \psi(a) \mathbf{k} &\text{ if }i=\bar{1}.\end{cases}\]
Obviously, $\psi_1$ becomes a group automorphism of $G$ of order $3$.
Now, we check that $\psi_1'$ becomes a group automorphism of $G'$ of order $3$. 
Note that the injectivity of $\psi_1'$ is straightforward. So $\psi_1'$ is bijective. 
For checking that $\psi_1'$ is a group homomorphism, we see that 
\[\psi_1' ((a_1,i_1)\cdot (a_2,i_2)) =\psi_1' (a_1\sigma_{i_1}(a_2) ,i_1+i_2)
=(\psi(a_1\sigma_{i_1}(a_2)\mathbf{i}^{i_1+i_2}) \mathbf{i}^{-(i_1+i_2)},i_1+i_2).\]
On the other hand, we have
\begin{align*}
\psi_1' (a_1,i_1)\cdot \psi_1'(a_2,i_2)
&=
(\psi(a_1 \mathbf{i}^{i_1}) \mathbf{i}^{-i_1},i_1)\cdot (\psi(a_2 \mathbf{i}^{i_2}) \mathbf{i}^{-i_2},i_2)\\
&=
(\psi(a_1 \mathbf{i}^{i_1}) \mathbf{i}^{-i_1} \sigma_{i_1} (\psi(a_2 \mathbf{i}^{i_2}) \mathbf{i}^{-i_2}),i_1+i_2)\\
&=
(\psi(a_1 \mathbf{i}^{i_1}) \mathbf{i}^{-i_1} \mathbf{i}^{i_1} \psi(a_2 \mathbf{i}^{i_2}) \mathbf{i}^{-i_2} \mathbf{i}^{-i_1},i_1+i_2)\\
&=
(\psi(a_1 \mathbf{i}^{i_1}) \psi(a_2 \mathbf{i}^{i_2}) \mathbf{i}^{-(i_1+i_2)},i_1+i_2)\\
&=
(\psi(a_1 \mathbf{i}^{i_1}a_2 \mathbf{i}^{i_2}) \mathbf{i}^{-(i_1+i_2)},i_1+i_2)\\
&=
(\psi(a_1 \mathbf{i}^{i_1}a_2 \mathbf{i}^{-i_1}\mathbf{i}^{i_1+i_2} ) \mathbf{i}^{-(i_1+i_2)},i_1+i_2)\\
&=
(\psi(a_1 \sigma_{i_1} (a_2)\mathbf{i}^{i_1+i_2}) \mathbf{i}^{-(i_1+i_2)},i_1+i_2).
\end{align*}
These imply $\psi_1' (a_1,i_1)\cdot \psi_1'(a_2,i_2)=\psi_1' ((a_1,i_1)\cdot (a_2,i_2))$.
We can also check that the order of $\psi_1'$ is $3$. In fact, 
\begin{align*}
\psi_1'^3 (a,i) =\psi_1' \circ \psi_1'(\psi(a \mathbf{i}^i) \mathbf{i}^{-i},i)
=\psi_1' (\psi^2(a \mathbf{i}^i) \mathbf{i}^{-i},i)
=(\psi^3(a \mathbf{i}^i) \mathbf{i}^{-i},i)= (a \mathbf{i}^i \mathbf{i}^{-i},i)=(a,i).
\end{align*}

Let $Q_1=(Q(G,\psi_1),*)$ and $Q_1'=(Q(G',\psi_1'),*')$ be the generalized Alexander quandles. 

For $Q_1$, since $(1,\bar{0})*(a,i)
=(a,i) \psi_1(a^{-1},i)
=(a,i) (\psi(a^{-1}),i)
=(a\psi(a^{-1}),\bar{0})$,
we have $P(Q_1)\subset Q_8 \times \{\bar{0}\}$.
Moreover, since we know $\mathbf{j}\psi(\mathbf{j}^{-1}) = -\mathbf{i}$ and $\mathbf{k}\psi(\mathbf{k}^{-1}) = -\mathbf{j}$
and $Q_8$ is generated by $-\mathbf{i}$ and $-\mathbf{j}$, we have $P(Q_1)=Q_8 \times \{\bar{0}\}$.
In a similar way, we can also obtain $P(Q_1')=Q_8 \times \{\bar{0}\}$.
Note that $P(Q_1')$ is not only $Q_8 \times \{\bar{0}\}$ as a set, but has the group structure $Q_8 \times \{\bar{0}\}$ as a direct product.
Thus both of the coset decompositions of $G/P(Q_1)$ and $G'/P(Q_1')$ are $Q_8 \times \{\bar{0}\} \sqcup Q_8 \times \{\bar{1}\}$.
Since neither of these two quandles $Q_1$ and $Q_1'$ satisfy (P2), 
we cannot use Theorem~\ref{thm:HK2} to determine whether $Q_1$ and $Q_1'$ are isomorphic. In the following, we use Theorem~\ref{thm:main} to show that $Q_1$ and $Q_1'$ are isomorphic.

Now we construct the two maps $h$ and $k$ in Theorem~\ref{thm:main} for $Q_1$ and $Q_1'$.
Since $P(Q_1)=P(Q_1')=Q_8 \times \{\bar{0}\}$ as sets,
we can define the map $h\colon P(Q_1) \to P(Q_1')$ to be the identity map.
Obviously, the condition (A) of Theorem~\ref{thm:main} can be verified.
We take $A=\{(1,\bar{0}),(1,\bar{1})\}$ and $A'=\{(1,\bar{0}),(-\mathbf{i},\bar{1})\}$
as the coset representatives of $G/P(Q_1)$ and $G'/P(Q_1')$,
respectively.
Then we define the map $k\colon A \to A'$ by
$k(1,\bar{0})=(1,\bar{0})$ and $k(1,\bar{1})=(-\mathbf{i},\bar{1})$.
Now we check the condition (B) of Theorem~\ref{thm:main}.
For $(1,\bar{0})\in A$, we have
\[
h((1,\bar{0})*(1,\bar{0}))=h(1,\bar{0})=(1,\bar{0})=(1,\bar{0})*'k((1,\bar{0})).
\]
For $(1,\bar{1})\in A$, we have
\[h((1,\bar{0})*(1,\bar{1}))=h((1,\bar{0}))=(1,\bar{0})\]
and
\begin{align*}
    (1,\bar{0})*'k(1,\bar{1})
    =(1,\bar{0})*' (-\mathbf{i},\bar{1})
    =(-\mathbf{i},\bar{1}) \psi_1' ((-\mathbf{i},\bar{1})^{-1})
    =(-\mathbf{i},\bar{1}) (\mathbf{i},\bar{1})
    =(-\mathbf{i}^2,\bar{0})=(1,\bar{0}).
\end{align*}
Hence the condition (B) of Theorem~\ref{thm:main} is also verified.
Therefore, by Theorem~\ref{thm:main}, we have $Q_1\cong Q_1'$.

Let $\tau \colon Q_8 \times \ZZ/2\ZZ \to Q_8 \times \ZZ/2\ZZ$ be the map defined by $\tau(a,i)=((-1)^i a ,i)$.
We can check that $\tau$ is 
both of group automorphisms of $G$ and $G'$.
Let $\psi_2= \tau \circ \psi_1$ and $\psi_2'= \tau \circ \psi_1'$.
Then both of the orders $\psi_2$ and $\psi_2'$ are $6$.

Let $Q_2$ and $Q_2'$ be the generalized Alexander quandles
$(Q(G,\psi_2),*)$ and $(Q(G',\psi_2'),*')$, respectively.

We construct $h$ and $k$ in Theorem~\ref{thm:main} to show 
$Q_2\cong Q_2'$. 

We can check $P(Q_2)=P(Q_2')=Q_8 \times \{\bar{0}\}$. 
Let $h: P(Q_2) \to P(Q_2')$ be the identity map on $Q_8 \times \{\bar{0}\}$.  The condition (A) of Theorem~\ref{thm:main} is satisfied. 
We take $A=\{(1,\bar{0}),(1,\bar{1})\}$ and $A'=\{(1,\bar{0}),(\mathbf{i},\bar{1})\}$
as the coset representatives of $G/P(Q_2)$ and $G'/P(Q_2')$,
respectively.
Then we define the map $k\colon A \to A'$ by
$k(1,\bar{0})=(1,\bar{0})$ and $k(1,\bar{1})=(\mathbf{i},\bar{1})$.
Now we verify the condition (B) of Theorem~\ref{thm:main}.
For $(1,\bar{0})\in A$, we have
\[
h((1,\bar{0})*(1,\bar{0}))=h(1,\bar{0})=(1,\bar{0})=(1,\bar{0})*'k((1,\bar{0})).
\]
For $(1,\bar{1})\in A$, we have
\[h((1,\bar{0})*(1,\bar{1}))
=h((1,\bar{1})\psi_2((1,\bar{1})^{-1}))
=h((1,\bar{1})\psi_2((1,\bar{1})))
=h((-1,\bar{1}+\bar{1}))
=(-1,\bar{0})
\]
and
\begin{align*}
    (1,\bar{0})*'k(1,\bar{1})
    =(1,\bar{0})*' (\mathbf{i},\bar{1})
    =(\mathbf{i},\bar{1}) \psi_2' ((\mathbf{i},\bar{1})^{-1}) 
    =(\mathbf{i},\bar{1}) (\mathbf{i},\bar{1})
    =(-1,\bar{0}).
\end{align*}
Hence the condition (B) of Theorem~\ref{thm:main} is satisfied. 
Therefore, by Theorem~\ref{thm:main}, we have $Q_2\cong Q_2'$.

\begin{rem}
For a quandle $Q= (Q, *)$ and a positive integer $n$, we can construct a quandle $Q^{(n)} = 
(Q, *^{(n)})$ by $ x *^{(n)} y = s_y^n (x)$ where $s_y$ is the point symmetry at $y$ for the quandle $(Q, *)$. 
If two qundles $Q$ and $Q'$ are isomorphic, then  $Q^{(n)}$ and $Q'^{(n)}$ are isomorphic (see \cite[Section 4]{HK2020}). 
For a generalized Alexander quandle $Q= Q(G, \psi)$, we have $Q^{(n)} = Q(G, \psi^n)$.  
For the quandles $Q_1$, $Q_1'$, $Q_2$ and $Q_2'$ above, we see that $Q_2^{(2)}\cong Q_1$ and $Q_2'^{(2)}\cong Q_1'$. 
Thus, there is no need to verify $Q_1\cong Q_1'$ once $Q_2\cong Q_2'$ is verified. However, for the sake of clarity of argument, we verified $Q_1\cong Q_1'$ at first. 
\end{rem}

\subsection{Complete list of $|\Qc_\gaq(n)|$ with $n \leq 127$}

Similar to the case of order $16$, there are some pairs of groups $G,G'$ and their automorphisms $\psi,\psi'$ 
such that we cannot use Theorem~\ref{thm:HK2} to determine whether $Q(G,\psi)$ and $Q(G',\psi')$ are isomorphic. 
(For example, we encounter such cases in order $24$ and $32$.) 
By using Theorem~\ref{thm:main}, we can classify $|\Qc_\gaq(n)|$ for small $n$'s. 
Table~\ref{fig:gaq} is the list of $|\Qc_\gaq(n)|$ with $n \leq 127$. 
These computations have been performed by using {\tt GAP} \cite{Gap} and the packages {\tt SmallGrp} \cite{SmallGrp} and {\tt Rig} \cite{Rig}. 
We have made code available at the following URL: 
\begin{center}
\url{https://github.com/Kurihara190/Classification_of_Generalized_Alexander_Quandles} 
\end{center}

\begin{rem}
The reason why the case $n = 128$ is not mentioned is that the number of nonisomorphic groups of order $128$ is significantly larger compared to orders up to $127$.
Consequently, given the computational capacity of our current computer, the search would require an impractical amount of time to complete the calculations. 
However, Theorem~\ref{thm:main} is valid for any order, so it may be possible to classify generalized Alexander quandles with orders higher than $128$ in the future, provided that computational capacity and efficiency are improved.
\end{rem}

\begin{table}[h]
    \centering
\begin{tabular}{|c|ccccccccccc|} \hline
    \textit{n} & 1 & 2 & 3 & 4 & 5 & 6 & 7 & 8 & 9 & 10 & 11 \\ 
    $|\mathcal{Q}_{\mathrm{GAQ}}(n)|$ & $\mathbf{1}$ & $\mathbf{1}$ & $\mathbf{2}$ & $\mathbf{3}$ & $\mathbf{4}$ & $\mathbf{3}$ & $\mathbf{6}$ & $\mathbf{9}$ & $\mathbf{11}$ & $\mathbf{5}$ & $\mathbf{10}$ \\ \hline 
    \textit{n} & 12 & 13 & 14 & 15 & 16 & 17 & 18 & 19 & 20 & 21 & 22 \\ 
    $|\mathcal{Q}_{\mathrm{GAQ}}(n)|$ & $\mathbf{11}$ & $\mathbf{12}$ & $\mathbf{7}$ & $\mathbf{8}$ & $\mathbf{29}$ & $\mathbf{16}$ & $\mathbf{17}$ & $\mathbf{18}$ & $\mathbf{15}$ & $\mathbf{13}$ & $\mathbf{11}$  \\ \hline
    \textit{n} & 23 & 24 & 25 & 26 & 27 & 28 & 29 & 30 & 31 & 32 & 33 \\ 
    $|\mathcal{Q}_{\mathrm{GAQ}}(n)|$ & $\mathbf{22}$ & $\mathbf{32}$ & $\mathbf{39}$ & $\mathbf{13}$& $\mathbf{51}$ & $\mathbf{20}$ & $\mathbf{28}$ & $\mathbf{15}$ & $\mathbf{30}$ & $\mathbf{87}$ & $\mathbf{20}$ \\ \hline
    \textit{n} & 34 & 35 & 36 & 37 & 38 & 39 & 40 & 41 & 42 & 43 & 44 \\ 
    $|\mathcal{Q}_{\mathrm{GAQ}}(n)|$ & $\mathbf{17}$ & $\mathbf{24}$ & $\mathbf{64}$ & $\mathbf{36}$ & $\mathbf{19}$ & $\mathbf{25}$ & $\mathbf{45}$ & $\mathbf{40}$ & $\mathbf{23}$ & $\mathbf{42}$ & $\mathbf{32}$ \\ \hline
    \textit{n} & 45 & 46 & 47 & 48 & 49 & 50 & 51 & 52 & 53 & 54 & 55 \\  
    $|\mathcal{Q}_{\mathrm{GAQ}}(n)|$ & $\mathbf{44}$ & $\mathbf{23}$ & $\mathbf{46}$ & $\mathbf{114}$ & $\mathbf{83}$ & $\mathbf{49}$ & $\mathbf{32}$ & $\mathbf{39}$ & $\mathbf{52}$ & $\mathbf{87}$ & $\mathbf{41}$ \\ \hline
    \textit{n} & 56 & 57 & 58 & 59 & 60 & 61 & 62 & 63 & 64 & 65 & 66 \\  
    $|\mathcal{Q}_{\mathrm{GAQ}}(n)|$ & $\mathbf{66}$ & $\mathbf{37}$ & $\mathbf{29}$ & $\mathbf{58}$ & $\mathbf{60}$ & $\mathbf{60}$ & $\mathbf{31}$ & $\mathbf{69}$ & $\mathbf{363}$ & $\mathbf{48}$ & $\mathbf{33}$ \\ \hline
    \textit{n} & 67 & 68 & 69 & 70 & 71 & 72 & 73 & 74 & 75 & 76 & 77 \\  
    $|\mathcal{Q}_{\mathrm{GAQ}}(n)|$ & $\mathbf{66}$ & $\mathbf{51}$ & $\mathbf{44}$ & $\mathbf{35}$ & $\mathbf{70}$ & $\mathbf{205}$ & $\mathbf{72}$ & $\mathbf{37}$ & $\mathbf{84}$ & $\mathbf{56}$ &$\mathbf{60}$ \\ \hline
    \textit{n} & 78 & 79 & 80 & 81 & 82 & 83 & 84 & 85 & 86 & 87 & 88 \\  
    $|\mathcal{Q}_{\mathrm{GAQ}}(n)|$ &$\mathbf{41}$ & $\mathbf{78}$ & $\mathbf{158}$ & $\mathbf{299}$ & $\mathbf{41}$ & $\mathbf{82}$ & $\mathbf{79}$ & $\mathbf{64}$ & $\mathbf{43}$ & $\mathbf{56}$ & $\mathbf{97}$ \\ \hline
    \textit{n} & 89 & 90 & 91 & 92 & 93 & 94 & 95 & 96 & 97 & 98 & 99 \\  
    $|\mathcal{Q}_{\mathrm{GAQ}}(n)|$ & $\mathbf{88}$ &$\mathbf{85}$ &$\mathbf{72}$ & $\mathbf{68}$ & $\mathbf{61}$ & $\mathbf{47}$ & $\mathbf{72}$ & $\mathbf{378}$ & $\mathbf{96}$ & $\mathbf{97}$ & $\mathbf{110}$ \\ \hline
    \textit{n} & 100 & 101 & 102 & 103 & 104 & 105 & 106 & 107 & 108 & 109 & 110 \\  
    $|\mathcal{Q}_{\mathrm{GAQ}}(n)|$ & $\mathbf{154}$ & $\mathbf{100}$ & $\mathbf{51}$ &$\mathbf{102}$ &$\mathbf{117}$ & $\mathbf{52}$ & $\mathbf{53}$ & $\mathbf{106}$ & $\mathbf{342}$ & $\mathbf{108}$ & $\mathbf{57}$ \\ \hline
    \textit{n} & 111 & 112 & 113 & 114 & 115 & 116 & 117 & 118 & 119 & 120 & 121 \\
    $|\mathcal{Q}_{\mathrm{GAQ}}(n)|$ & $\mathbf{73}$ & $\mathbf{204}$ & $\mathbf{112}$ & $\mathbf{59}$ & $\mathbf{88}$ &$\mathbf{87}$ &$\mathbf{135}$ & $\mathbf{59}$ & $\mathbf{96}$ & $\mathbf{177}$ & $\mathbf{219}$ \\ \hline
    \textit{n} & 122 & 123 & 124 & 125 & 126 & 127 &  &  & & & \\ 
    $|\mathcal{Q}_{\mathrm{GAQ}}(n)|$ & $\mathbf{61}$ & $\mathbf{80}$ & $\mathbf{92}$ & $\mathbf{297}$ & $\mathbf{127}$ & $\mathbf{126}$ & & & & & \\ \hline 
\end{tabular} \\
    \caption{List of $|\Qc_\gaq(n)|$ with $n \leq 127$}\label{fig:gaq}
\end{table}

\bigskip

\noindent 
\textbf{Author's addresses.}

\smallskip

\noindent
Akihiro Higashitani: \\
Department of Pure and Applied Mathematics, 
Osaka University, Suita, Osaka 565-0871, Japan\\
E-mail address: {\tt higashitani@ist.osaka-u.ac.jp}

\smallskip

\noindent
Seiichi Kamada: \\
Department of Mathematics, 
Osaka University, Toyonaka, Osaka 560-0043, Japan\\
E-mail address:  {\tt kamada@math.sci.osaka-u.ac.jp}

\smallskip

\noindent
Jin Kosaka: \\
Department of Mathematics, 
Osaka University, Toyonaka, Osaka 560-0043, Japan\\
E-mail address:  {\tt kosakajin123@gmail.com}

\smallskip

\noindent
Hirotake Kurihara: \\
Department of Applied Science, 
Yamaguchi University, Ube 755-8611, Japan\\ 
E-mail address: {\tt kurihara-hiro@yamaguchi-u.ac.jp}

\end{document}